\documentclass[review,10]{elsarticle}

\usepackage{lipsum}
\usepackage{hyperref}
\usepackage[hyperref]{xcolor}

\hypersetup{
  colorlinks=true,
}
\usepackage[center]{caption}
\usepackage{xcolor}
\usepackage{amsmath}
\usepackage{amsfonts}
\usepackage{graphicx}
\usepackage{subfigure}
\usepackage{booktabs}
\usepackage{empheq}
\usepackage[
singlelinecheck=false 
]{caption}
\newtheorem{thm}{Theorem}[section]
\newtheorem{lem}{Lemma}[section]
\newdefinition{rmk}{Remark}[section]
\newdefinition{remark}{Remark}[section]

\newproof{pf}{Proof}
\renewenvironment{pf}{\noindent{\textit{Proof.}  }}{\hfill $\Box$}

\def\RR{\mathbb{R} }

\def\EE{\mathbb{E}}

\def\Cov{{\rm Cov}}

\newproof{pot1}{Proof of Theorem \ref{th:1}}
\newproof{pfL1}{Proof of Lemma \ref{e:Lop}}
\newproof{pfL2}{Proof of Lemma \ref{lem:1}}
\newproof{pfL3}{Proof of Lemma \ref{lem:2}}
\numberwithin{equation}{section}

\begin{document}

\begin{frontmatter}

\title{Estimation of all parameters in the reflected Orntein-Uhlenbeck process from discrete observations}


\author[mymainaddress]{Yaozhong Hu}
\ead{yaozhong@ualberta.ca}

\author[mysecondaryaddress]{Yuejuan Xi\corref{mycorrespondingauthor}}
\cortext[mycorrespondingauthor]{Corresponding author}
\ead{yjx@mail.nankai.edu.cn}

\address[mymainaddress]{Department of Mathematical and Statistical Sciences, University of Alberta at Edmonton Edmonton, Alberta Canada, T6G 2G1}

\address[mysecondaryaddress]{School of Mathematical Sciences, Nankai University, Tianjin, PR China, 300071}

\begin{abstract}
Assuming that a reflected Ornstein-Uhlenbeck state process 
is observed at discrete time instants, we propose  generalized  moment estimators to estimate all drift and diffusion parameters via the celebrated ergodic theorem. With the sampling time step $h>0$ arbitrarily fixed, we prove the strong consistency and asymptotic normality of our  estimators as the sampling  size $n$ tends to infinity. 
This   provides   a complete solution to an open problem left in 
\citet{MR3395608}. 
\end{abstract}


\end{frontmatter}

\section{Introduction}
On  a filtered probability space $ (\Omega, \mathbb P, \mathcal F, 
\{\mathcal F_t\}_{t\ge 0} )$  let $W=\{W(t)\}_{t\ge0}$ be a one-dimensional standard Brownian motion. All the processes mentioned in this paper will be adapted to  $\{\mathcal F_t\}_{t\ge 0} $.  We consider the (reflected) Ornstein-Uhlenbeck  (ROU) process, reflected at zero, which is defined by the following one-dimensional stochastic differential equation (SDE):
\begin{equation}\label{e:ROU}
\begin{cases}
dX_t=\kappa (\theta-X_t)dt+\sigma dW_t+L_t,\quad t\in \RR_+=\left\{x, x\ge 0\right\} \,, \\
X_0=x \in \RR_+\,, \\
\end{cases} 
\end{equation}
where $\kappa, \theta, \sigma\in(0,\infty)$ are  constants  and 
 $L_t$ is the minimal continuous increasing process which ensures that $X_t\geq0$ for all $t\geq0$. The ROU process is a useful stochastic model in  finance and queue theory  (cf. \citet{MR2144561},  \citet{MR1957808} and the references therein). 

This paper will concern with  the statistical  estimation problem for the 
parameters $\kappa, \theta, \sigma$ from the observations.
In most practical situations  the observations of the process $\{X_t, t\ge 0\}$
can be made only at discrete 
  time  instants $t_k=kh$, $k=1, 2, \cdots, $ and   usually the time interval 
  $h$ between consecutive observations cannot be made arbitrarily small. 
  To deal with this situation an   ergodic type of estimator to estimate 
  $\kappa$ and $\theta$ is proposed in 
  a previous work \citet{MR3395608} and the strong consistency and asymptotic normality of the estimators are also obtained there. 
  However, as pointed out in \citet{MR3395608}  they were unable to estimate
  $\sigma$ (or $\sigma^2$) by using the ergodic type estimator 
  and instead they  proposed to use 
  $\hat \sigma_{c, n} :=\frac{1}{nh}\sum_{k=1}^n(X_{(k+1)h}-X_{kh})^2$ as the estimator of $\sigma^2$.
Let us also mention a work on the estimation of the parameters  $\kappa, \theta, \sigma$ for  this  ROU when continuous observation is available (\citet{MR2719520}). 
This would require that $h\rightarrow 0$ to guarantee   the 
strong consistency  of the estimator
(e.g. $\hat \sigma_{c, n}^2 \rightarrow \sigma^2$).  This paper will fill this gap.
We shall introduce an ergodic type estimator to estimate $\sigma^2$
(and hence we can estimate all the parameters $\kappa, \theta, \sigma^2$  simultaneously) and prove the strong consistency and asymptotic 
normality for all estimators (including the estimator $\hat \sigma_n$ for
$\sigma$) regardless the (fixed) value of $h$.  This work is motivated by a recent work of 
  \citet{MR4072252},  where the Ornstein-Uhlenbeck process has no reflection, but the Brownian motion was replaced by a stable process. 
   
Now let us   describe  our ergodic estimators   for all parameters. 
 It is well-known that 
there is a unique invariant probability density function $\pi(x)$ of $X_t$ such that for any  integrable function $f$ 
we have 
\begin{equation}\label{e:Inv1}
\lim_{n\rightarrow \infty}\frac{1}{n} \sum_{k=1}^n f(X_{t_k})=\int_{\RR_+} f(x) \pi(x) dx ,
\end{equation}
and the invariant probability density function $\pi(x)$ has the following explicit expression (see \citet{MR3395608}): 
\begin{equation}\label{e:sta}
\pi(x)=\frac{\sqrt{2\kappa}}{\sigma}\frac{\phi(\frac{\sqrt{2\kappa}(x-\theta)}{\sigma})}{1-\Phi(-\frac{\sqrt{2\kappa}\theta}{\sigma})},
\end{equation}
 $\phi(x)=e^{-\frac{|x|^2}{2} }/{\sqrt{2\pi}}$ is the standard normal probability density function, and $\Phi(x)=\int_{-\infty}^x \phi(u) du$ is the standard normal distribution function. As observed in \citet{MR3395608}   
 the invariant measure $\pi(x)$ remains  the same function if  
 the quantities $\theta$ and   $\frac{ \kappa }{\sigma^2}$ remain unchanged. 
 Thus, we cannot expect to use \eqref{e:Inv1} to estimate $\kappa$ and $\sigma^2$ simultaneously. 
To this end and motivated by  \citet{MR4072252} we  shall  use the ergodic theorem for $ X_{kh}, X_{(k+1)h}$,  which states that for any integrable function $f:\RR_+^2\rightarrow \RR$,
\begin{equation}
\lim_{n\to\infty}\frac{1}{n}\sum_{k=1}^nf(X_{kh}, X_{(k+1)h})=\EE f(\tilde X_0, \tilde X_h)=\int_{\RR_+^2} f(x,y)\pi(x)p_h(x,y)dxdy,
\label{e.1.4}
\end{equation}
where $\tilde X_0$ is a random variable independent of the Brownian motion $W$
and having the invariant probability density $\pi(x)$, $\tilde X $ is the solution to 
\eqref{e:ROU} with initial random variable $\tilde X_0$, and $p_h(x,y)$ is the transition density of $X$. 

With some specific choices of $f$ in \eqref{e:Inv1}  and \eqref{e.1.4} we can
obtain our ergodic estimators, whose detailed construction  is given in the
next section, where the strong consistency and asymptotic normality are also obtained.

Section \ref{sec:Num} will provide a numerical example which demonstrates the convergence 
results of our estimators and which also demonstrates that $\hat \sigma_{c, n}$ 
does not converge.


\section{Strong consistency and asymptotic normality}\label{sec:2}
In this section, we aim to construct the estimators  for all the parameters $\kappa, \theta, \sigma$ of the ROU process $\{X_t, t\ge 0\}$ 
given by \eqref{e:ROU} 
  based on discrete observations $\{X_{t_1,}, \cdots, X_{t_n}\}$,
  where $t_k=kh$ with the observation time interval $h$ arbitrarily fixed.
We will also  study their strong consistency and asymptotic normality. 
We  begin with two crucial convergence results, which are adapted from Lemma 1 in \citet{MR3395608} and Theorem 1.1 in \citet{billingsley1961statistical}, respectively.
\begin{lem}\label{le:1}
The $h$-skeleton sampled chain $\{X_{kh}:k\ge0\}$ is ergodic. Namely, for any 
initial value $x\in\mathcal \RR_+$ and $f\in L_1(\RR_+)$
and $g\in L_1(\RR_+^2)$
\ we have 
\begin{empheq}[left=\empheqlbrace]{align} 
&\lim_{N\to\infty}\frac{1}{n}\sum_{k=1}^n f(X_{kh})= \mathbb E[f(X_\infty)]=
\int_{\mathbb R_+} f(x)\pi(x)dx, \quad a.s.,\label{e.2.1} \\
&\lim_{n\to\infty}\frac{1}{n}\sum_{k=1}^n g(X_{kh}, X_{(k+1)h}) =\EE g(\tilde X_0, \tilde X_h)\nonumber\\
&\qquad \qquad \qquad \qquad =\int_{\RR_+^2  } g(x,y)\pi(x)p_h(x,y)dxdy,\; \quad a.s., \label{e.2.2} 
\end{empheq}
where 
$\tilde X_0$ is a random variable independent of the Brownian motion $W$
and having the invariant probability density $\pi $, $\tilde X $ is the solution to 
\eqref{e:ROU} with initial random variable $\tilde X_0$, and 
$p_h(x,y)$ is the transition density of $X$.
\end{lem}
As illustrated in  \citet{MR3395608} it is impossible to use 
\eqref{e.2.1} alone to estimate all the parameters $\kappa, \theta, \sigma$.
So   we take $f_1(x)=x$ and $f_2(x)=x^2$ in \eqref{e.2.1}  and we take 
$g(x,y)=xy$ in \eqref{e.2.2} to obtain a system of three equations to determine the parameters $\kappa, \theta, \sigma$.
Some elementary computations yield the following expressions for 
 the stationary moments of the invariant measure.
\begin{equation}\label{e:2}
    \left\{
     \begin{aligned}
        &\EE [X_\infty]= \theta + \frac{\sigma}{\sqrt{2\kappa}} \frac{\phi(\frac{\sqrt{2\kappa }\theta}{\sigma})}{1-\Phi(-\frac{\sqrt{2\kappa }\theta}{\sigma})}\,,\\
        &\EE [X_\infty^2]= \frac{\sigma^2}{2\kappa}+\theta^2+\theta\frac{\sigma}{\sqrt{2\kappa}} \frac{\phi(\frac{\sqrt{2\kappa }\theta}{\sigma})}{1-\Phi(-\frac{\sqrt{2\kappa }\theta}{\sigma})}\,, \\
       & \EE   (\tilde X_0\tilde X_h)=
\int_{\mathbb R_+^2}  xy \pi(x)p_h(x,y)  dx dy\,. 
     \end{aligned}
     \right.
\end{equation}
However, to our best knowledge, there is no compact explicit form for the transition 
probability density $p_t(x,y)$.  
We  shall use the following 
spectral representation for  the transition density 
$p_t(x,y)$  derived in \citet{MR2144561} 
\begin{equation*}
p_t(x,y)=\pi(y)+m(y)\sum_{i=1}^\infty e^{-\lambda_i t}\varphi_i(x)\varphi_i(y), \quad t>0,
\end{equation*}
where  the notations are described as follows.
\begin{enumerate}
\item[(1)] $\pi(x)$ is   the stationary density given by \eqref{e:sta}  and 
 $m(x)$ is the speed measure   defined by
\begin{equation*}
m(x)=\frac{2}{\sigma^2 }e^{-\kappa(\theta-x)^2/\sigma^2}\,. 
\end{equation*}
\item[(2)] The eigenvalues 
$0<\lambda_1<\lambda_2<\cdots<\lambda_i<\cdots $ are roots of 
\begin{equation*}
H_{\lambda/\kappa-1}(-\sqrt\kappa\theta/\sigma)=0\,, 
\end{equation*}
where $H$ is the Hermite function   (see \citet{MR0174795}). 
\item[(3)] The normalized eigenfunctions $\varphi_i(x)$ are given by
\begin{equation*}
\varphi_i(x)=\pm \frac{\kappa^{3/4}\sigma^{1/2} e^{\kappa\theta^2/(2\sigma^2)}H_{\lambda_i/\kappa}(\sqrt\kappa(x-\theta)/\sigma)}{\sqrt{2\lambda_i \bigtriangleup_iH_{\lambda_i/\kappa}(-\sqrt\kappa\theta/\sigma)}},\; i=1,2,\cdots
\end{equation*}
where $\bigtriangleup_i=\frac{\partial H_{\nu-1}(-\sqrt\kappa\theta/\sigma)}{\partial \nu}|_{\nu=\lambda_i/\kappa}$.
\end{enumerate} 
Now we replace $\EE(X_\infty), \EE(X_\infty^2), \EE(\tilde X_0\tilde X_h)$
in \eqref{e:2} by their  sample approximations to yield 
\begin{equation}\label{e.2.4}
    \left\{
     \begin{aligned}
        &\frac{1}{n} \sum_{k=1}^n {X_{kh}}= \theta + \frac{\sigma}{\sqrt{2\kappa}} \frac{\phi(\frac{\sqrt{2\kappa }\theta}{\sigma})}{1-\Phi(-\frac{\sqrt{2\kappa }\theta}{\sigma})},\\
        &\frac{1}{n} \sum_{k=1}^n  X_{kh}^2= \frac{\sigma^2}{2\kappa}+\theta^2+\theta\frac{\sigma}{\sqrt{2\kappa}} \frac{\phi(\frac{\sqrt{2\kappa }\theta}{\sigma})}{1-\Phi(-\frac{\sqrt{2\kappa }\theta}{\sigma})},\\
       & \frac{1}{n} \sum_{k=1}^n X_{kh}X_{(k+1)h}=
\int_{\mathbb R_+^2}  xy \pi(x)p_h(x,y)  dx dy.
     \end{aligned}
     \right.
\end{equation}
This is a system of three equations for the three unknown parameters. 
We expect that it would give a unique solution 
$\hat  \kappa_n, \hat\theta_n, \hat\sigma_n$, 
which we call the ergodic estimators
of the parameters. The system is still complicated to analyze and to be solved.
We will further simplify it.
To this end we  denote $u=\theta$ and $v=\frac{\sqrt{2\kappa}\theta}{\sigma}$. Then the first two equations in \eqref{e.2.4} depends only on
$u$ and $v$  and  they give a unique solution 
$\hat u_n$ and $\hat v_n$.  We then write  $p_h(x,y)=p_h(x,y; \kappa, \theta, \sigma)$ as 
a kernel $p_h(x,y; u, v, \sigma)$ depending on parameters $u,v, \sigma$.
Finally, we replace the parameters $u$ and $ v$ in 
kernel $p_h(x,y; u, v, \sigma)$   by the 
obtained  values $\hat u_n$ and $\hat v_n$,  then the third equation in \eqref{e.2.4}
becomes one equation for one unknown $\sigma$.
This greatly simplifies the computations.

To summarize the above discussion, we have transformed the system \eqref{e.2.4} into
the following system of equations.
\begin{equation}\label{e:system3}
    \left\{
     \begin{aligned}
        &\frac{1}{n}\sum_{k=1}^n X_{kh}= u + \frac{u}{v} \frac{\phi(v)}{1-\Phi(-v)},\\
        &\frac{1}{n}\sum_{k=1}^n X_{kh}^2= \frac{u^2}{v^2}+u^2+\frac{u^2}{v}\frac{\phi(v)}{1-\Phi(-v)},\\
        &\frac{1}{n}\sum_{k=1}^n X_{kh}X_{(k+1)h}= \int_0^\infty\int_0^\infty xy  p_h(x,y) \pi(x) dxdy\,.
     \end{aligned}
     \right.
\end{equation}
The  right-hand side of the above 
third equation depends on $u$, $v$ and $\sigma$ by substituting 
$\theta$ and $\kappa$  by $\theta=u$ and $\kappa=\frac{v^2\sigma^2}{2u^2}$ into the expression of $\pi$ and
$p_h(x,y)$. Define $\tilde\lambda_i=\lambda_i/\sigma^2$. For sake of the numerical computation 
we write the dependence explicitly as follows:
\begin{equation}
 p_h(x,y)= \pi(y)+ m(y)\sum_{i=1}^\infty e^{-\tilde\lambda_i \sigma^2 h}\varphi_i(x) \varphi_i(y),\label{e.2.6}
\end{equation}
where 
\begin{enumerate}
 \item[(1)]\   $ \pi$ and $  m$ are  given by 
 \begin{equation}
 \pi(x)=\frac{ v  }{ u  }\phi\left(\frac{ v  }{ u  }x- v  \right)/[1-\Phi(- v  )],\label{e.2.7}
\end{equation}
and 
\begin{equation}
 m(x)=\frac{2}{\sigma^2}e^{- v  ^2/2+ v  ^2x/ u  - v  ^2x^2/(2 u  ^2)},\label{e.2.8}
\end{equation}
\item[(2)] \ The eigenvalues $0<\tilde\lambda_1<\tilde\lambda_2<\cdots<\tilde\lambda_i<\cdots $ are roots of 
\begin{equation}
H_{2 u  ^2\tilde\lambda/ v  ^2-1}(- v  /\sqrt 2)=0\,.
\label{e.2.9}
\end{equation}  
\item[(3)] The eigenfunctions are given by 
\begin{equation}
\varphi_i(x)= \pm\sigma\frac{( v  /(\sqrt 2  u  ))^{3/2}e^{ v  ^2/4}H_{2 u  ^2\tilde\lambda_i/ v  ^2}((\frac{ v  }{ u  }x- v  )/\sqrt 2)}{\sqrt{2\tilde \lambda_i  \bigtriangleup_i H_{2 u  ^2\tilde\lambda_i/ v  ^2}(- v  /\sqrt 2)}} ,\; \quad i=1,2,\cdots
\label{e.2.10}
\end{equation}
with $ \bigtriangleup_i=\frac{\partial H_{\nu-1}(- v  /\sqrt 2)}{\partial \nu}|_{\nu=2 u  ^2\tilde\lambda_i/ v  ^2}$
\,.  
\end{enumerate}
 Now we summarize  our discussion as follows. 

\noindent{\bf Construction of the ergodic estimators for all parameters 
$\kappa, \theta, \sigma$}:
\begin{enumerate}
\item[(i)]  Solve   the first  two equations in the system 
\eqref{e:system3} to obtain $\hat u_n$ and $\hat v_n$.  
\item[(ii)] Substitute the obtained  $\hat u_n$ and $\hat v_n$  into the transition probability kernel $p_h(x,y)$ according to 
\eqref{e.2.6}-\eqref{e.2.10}   to obtain the third
equation  in the system  \eqref{e:system3}, 
which now 
contains only one unknown $\sigma$  and 
  solve it  to obtain $\hat \sigma_n$.   
\item[(iii)] Solve $\hat u_n=\hat \theta_n  \quad{\rm and}\quad 
\hat v_n=\frac{\sqrt{2\hat\kappa_n }\hat\theta_n}{\hat\sigma_n}$ to obtain 
\begin{equation} 
 \hat \theta_n=\hat u_n  \quad{\rm and}\quad 
\hat \kappa_n= \frac{\hat v_n^2 \hat \sigma_n^2}{
2\theta_n^2}  \,.\label{e.2.11} 
\end{equation} 
\end{enumerate} 
\begin{remark}
In numerical computation, we shall need to take finite terms in the spectral 
representation of the transition probability function
(in our numerical simulation we take about  twelve terms and the results are satisfactory). The  Hermite 
functions and the roots of the Hermite functions can be handled by the 
standard mathematical software package. The system \eqref{e:system3} of algebraic equations does not give an explicit solution.  There are many standard methods to solve it, such as the Newton-Raphson iteration method. 
\end{remark}
To study the strong consistency and the asymptotic normality,
we  denote the right-hand sides of the equation in the system 
\eqref{e:system3} by $g_1(u,v)$, $g_2(u,v)$, and $g_3(u,v,\sigma)$, respectively. Denote  
\begin{equation*}
M_{1,n}=\frac{1}{n}\sum_{k=1}^nX_{kh},  \quad 
 M_{2,n}=\frac{1}{n}\sum_{k=1}^n(X_{kh})^2, \; \quad  M_{3,n}=\frac{1}{n}\sum_{k=1}^n X_{kh}X_{(k+1)h}\;.
\end{equation*}
Then the equation \eqref{e:system3} can be rewritten as 
\begin{equation}
g_1(u,v)=M_{1,n}, \quad g_2(u,v)=M_{2,n}, \quad g_3(u,v, \sigma)=M_{3,n}\,. 
\label{e.2.12}  
\end{equation}
($g_3$ also depends on $h$ which is fixed) Or we write 
\begin{equation}
g(u,v,\sigma)=M_n\,,  \label{e.2.13}
\end{equation}
where
\[ 
g=(g_1, g_2, g_3)^T\quad{\rm and} \quad M_n=(M_{1,n}, M_{2,n},M_{3,n})^T\,.
\]
Denote by $J(u,v, \sigma)$ the determinant of the Jacobian of $g$.  Then 
\begin{equation}
J(u,v, \sigma)=J(g_1, g_2) \frac{\partial}{\partial \sigma}g_3
( u ,   v, \sigma)\,, 
\end{equation} 
where $J(g_1, g_2)$ is the determinant of the Jacobian of $g_1$ ad $g_2$.
\citet{MR3395608} proved that $J(g_1, g_2)$ is never $0$.
If $\frac{\partial}{\partial \sigma}g_3
( u  ,   v , \sigma)$ is not singular in some domain $D\subseteq \RR_+^3$,
 then by the inverse  function 
theorem, for any $(u,v, \sigma)\in D$,
 $g=(g_1, g_2, g_3)$ has a   unique inverse in a neighbourhood $(u, v, \sigma)$.
If $(u, v, \sigma)$  (or equivalently, $(\kappa, \theta, \sigma)$) are the true parameters,  then by Lemma \ref{le:1}, we see when $n$ is sufficiently 
large $( M_{1,n}, M_{2,n}, M_{3,n})$ 
will be in the  neighbourhood of $g(u,v, \sigma)$.  This means when $n$ is sufficiently large the equation \eqref{e.2.12} has a solution. 

Thus,  the critical question now  is to find a domain $D$ such that 
$\frac{\partial}{\partial \sigma}g_3
( u  ,   v , \sigma)$ is not singular on $D$. 
This is an elementary analysis problem.
  The explicit expression of the derivative 
of $g_3(u,v, \sigma)$ with respect to $\sigma$ can be obtained
(see Remark \ref{r.2.1} below). However, this expression is complicated and 
 it is hard to obtain the  domain of $(u, v, \sigma)$ 
so that inside this domain    this derivative is not  singular.
We shall proceed as follows to reduce the 
$\frac{\partial}{\partial \sigma}g_3
( u  ,   v , \sigma)$ from a function of three variables 
$u  ,   v , \sigma $  to a function of one variable $\sigma$.

Since the first two equations in \eqref{e:system3} is independent of $\sigma$, 
as indicated above we can solve them without considering the third equation
in \eqref{e:system3}. 
 \citet{MR3395608} proved that there exist continuous inverse  mapping
  $
 (h_1, h_2)$ of $(g_1, g_2):
 \RR_+^2\rightarrow \RR^2$  such that  the ergodic estimators defined by
\begin{equation}\label{e:uv}
\hat u_n:=h_1(M_{1,n},M_{2,n}),\qquad   \hat v_n:=h_2(M_{1,n},M_{2,n})
\end{equation}
converge almost surely to the true parameters 
\begin{equation*}
u=h_1(g_1(u,v),g_2(u,v))=\theta ,\qquad  v=h_2(g_1(u,v),g_2(u,v))=\frac{\sqrt{2\kappa}\theta}{\sigma}\,. 
\end{equation*}
After  the estimators $\hat u_n$ and $\hat v_n$
have been obtained, we can substitute them into the $g_3 $.  Thus $g_3
(\sigma)=g_3
(\hat u_n, \hat v_n, \sigma)$ and $g_3'
(\sigma)=\frac{\partial}{\partial \sigma}g_3
(\hat u_n, \hat v_n, \sigma)$  will be   functions of single variable
$\sigma$. We can plot the  derivative function $g_3'
(\sigma)$ in an interval  $D_\sigma$ that is  as large as we believe it contains   
the true parameter $\sigma$
(we shall plot $g_3'
(\sigma)$ for some value of $u$ and $v$ in next section). If   
$g_3'$ is never equal  to $0$ on $D_\sigma$,  then the solution 
to  $g_3(\hat u_n, \hat v_n, 
\sigma)=M_{3, n}$    is   unique on $D_\sigma$   (if not then there are two different points 
$\sigma_1<\sigma_2$ in $D_\sigma$ such that $g_3(\hat u_n, \hat v_n, 
\sigma_1)=g_3(\hat u_n, \hat v_n, 
\sigma_2)=M_{3, n}  $. By the mean value theorem there is a $\sigma_0\in 
[\sigma_1, \sigma_2]\subseteq 
D_\sigma$ such that $\frac{\partial}{\partial \sigma}g_3
(\hat u_n, \hat v_n, \sigma)=0$).

If   $g_3'$ is not singular  on $D_\sigma$, then
the third equation  \eqref{e.2.12} has a unique solution 
on $D_\sigma$, which gives the ergodic estimator 
$\hat \sigma_n=h_3(\hat u_n, \hat v_n, M_{3,n})$ of $\sigma$,
where $h_3(\hat u_n, \hat v_n, \cdot)$ is the continuous inverse of $
g (\hat u_n, \hat v_n, \cdot)$.  By Lemma \ref{le:1} 
it is easy to see that $\hat \sigma_n \rightarrow \sigma$ a.s.  

\begin{figure}
\centering 
\includegraphics[width=0.8\textwidth]{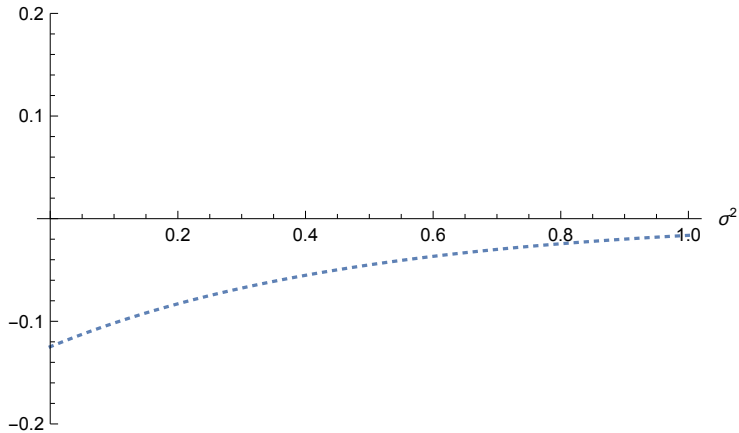} 
\caption{
The graph of $\frac{1}{h}\frac{\partial }{\partial \sigma^2}g_3(u,v,\sigma)$ 
for $u=1$, $v=\sqrt2/0.5$, $h=0.5$. 
} 
\label{Fig.1} 
\end{figure}

  Now we summarize the above discussion as the following theorem.
\begin{thm}
\begin{enumerate}
\item[(i)]  The first two equations of the system \eqref{e:system3} have a 
unique solution pair $(\hat u_n, \hat v_n)=(h_1(M_{1,n}, M_{2, n}), 
h_2(M_{1,n}, M_{2, n}))$. 
\item[(ii)]  If $\frac{\partial}{\partial \sigma}g_3
(\hat u_n, \hat v_n, \sigma ) $     ($g_ 3$ is defined by 
the  right-hand side of the third equation in \eqref{e:system3})
is not singular on some interval $\sigma\in D_\sigma$ which contains the true parameter 
$\sigma$, then when $n$ is sufficiently large 
the third equation of \eqref{e:system3}, namely, 
\begin{equation}
g_3(\hat u_n, \hat v_n, \sigma)=M_{3,n}=\sum_{k=1}^n X_{k h} X_{(k+1) h} 
\end{equation}
has a unique solution $\hat \sigma_n$. 
\item[(iii)] $ (\hat \kappa_n, \hat \theta_n, \hat
\sigma_n )^T\rightarrow (\kappa, \theta, \sigma)^T$ almost surely as $n\rightarrow \infty$,
where  $ 
 \hat \theta_n$ and $
\hat \kappa_n$ are given by \eqref{e.2.11}. 
\end{enumerate} 
\end{thm}
Next, we study the joint asymptotic behavior of the all estimators $(\hat \theta_n, \hat \kappa_n, \hat \sigma_n)$.
\begin{thm}
Let  $\frac{\partial}{\partial \sigma}g_3
(\hat u_n, \hat v_n, \sigma ) $   be  nonsingular on some interval $\sigma\in D_\sigma$ which contains the true parameter 
$\sigma$. 
 Then,  the estimators  $(\hat\theta_n, \hat\kappa_n,\hat\sigma_n)$ satisfy the following asymptotic normality property: 
\begin{equation*}
\sqrt n((\hat\theta_n,\hat\kappa_n,\hat\sigma_n)^T-(\theta,\kappa,\sigma)^T)\Rightarrow N(0,\Sigma),
\end{equation*}
where $\Sigma$ is a covariance matrix defined in \eqref{e:Sig} below.
\end{thm}
\begin{pf}
For any nice function $f$ and $g$, denote
\begin{align*}
\sigma_{fg}:=&\Cov (f(\tilde X_0,\tilde X_h), g(\tilde X_0,\tilde X_h))+\sum_{k=1}^\infty [\Cov (f(\tilde X_0, \tilde X_h), g(\tilde X_{kh}, \tilde X_{(K+1)h}))\\
&+\Cov (g(\tilde X_0, \tilde X_h), f(\tilde X_{kh}, \tilde X_{(k+1)h}))].
\end{align*}
Let $f_1(x,y)=x$, $f_2(x,y)=y$, $f_3(x,y)=xy$ 
and denote 
\begin{equation*}
\tilde\Sigma_3:=(\sigma_{f_kf_l})_{1\leq k,l\leq 3}.
\end{equation*}
Then an application of  the multivariate Markov chain 
central limit theorem  (e.g.    \citet[Section 1.8.1]{MR2742422})  yields 
\begin{equation*}
\sqrt n((M_{1,n}, M_{2,n}, \tilde M_{2,n})^T-(\EE X_\infty, \EE X_\infty ^2, \EE \tilde X_0 \tilde X_h )^T) \stackrel{d}\rightarrow 
N(0,\tilde \Sigma_3).
\end{equation*}
To simplify  notations, introduce the following  two mappings:
\begin{equation*}
h:(x_1,x_2,x_3)\mapsto (h_1(x_1,x_2), h_2(x_1,x_2), h_3(x_1,x_2,x_3)), 
\end{equation*}
and
\begin{equation*}
\eta:(x_1,x_2,x_3)\mapsto  (x_1, \frac{x_2^2x_3^2}{2x_1^2}, x_3),
\end{equation*}
where $\eta$ is the inverse transform of \eqref{e.2.11}. Then from    delta method, we have
\begin{equation*}
\sqrt n( h(M_{1,n}, M_{2,n}, \tilde M_{2,n})^T-h(\kappa, \theta, \sigma) )^T  \stackrel{d}\rightarrow N(0, \bar\Sigma)
\end{equation*}
where $ \bar\Sigma=\nabla h(\Theta)\tilde\Sigma_3 \nabla h(\Theta)^T$.
Finally, applying the  delta method again, we arrive at the asymptotic 
behavior of the ergodic  estimators: 
\begin{equation*}
\sqrt n(\eta( h(M_{1,n}, M_{2,n}, \tilde M_{2,n}))^T-\eta (h(\Theta))^T )\Rightarrow N(0, \Sigma)
\end{equation*}
where
\begin{equation}\label{e:Sig}
 \Sigma=\nabla \eta(h(\Theta))\tilde\Sigma_3\nabla \eta(h(\Theta))^T 
 \end{equation}
  completing  the proof of the theorem.
\end{pf}
\begin{rmk}\label{r.2.1}
 We mention that we do not know the monotonicity of $g_3^\prime(\sigma)$ in theory. But we can observe it numerically. Since $\sigma>0$, to investigate the sign of $g_3^\prime(\sigma)$ is equivalent to discuss the sign of $\frac{\partial }{\partial \sigma^2}g_3(u,v,\sigma)$. 
 \begin{equation*}
 \frac{\partial }{\partial \sigma^2}g_3(u,v,\sigma)=-h\int_0^\infty \int_0^\infty xy[m(y)\sum_{i=1}^\infty \tilde\lambda_ie^{-\tilde
 \lambda_i\sigma^2h}\varphi_i(x)\varphi_i(y)]dxdy.
 \end{equation*}
  For fixed $u$ and $v$ an example of the values of $\frac{1}{h}\frac{\partial }{\partial \sigma^2}g_3(u,v,\sigma)$  is plotted  in Figure \ref{Fig.1}. It shows that the partial derivatives are always less than zero on the concerned interval.  
\end{rmk}
\section{Numerical experiments}\label{sec:Num}
In this section, we present a numerical experiment to illustrate our method. 
In Table \ref{tab:1} we set the following true parameters $
\sigma=0.5$, $\kappa=1$, $\theta=1$. The time step is fixed by $h=0.5$. In the experiments, we use the
truncation
\begin{equation*}
 p_{N,h}(x,y)=\pi(x)+\sum_{i=1}^{N} e^{-\tilde\lambda_i \sigma^2 h}\varphi_i(x)\varphi_i(y)
\end{equation*}
in \eqref{e.2.6}. Here we take $N=12$. It can be seen that our estimators for all the parameters, including
$\hat\sigma_n$ are strongly consistent. 
On the other hand  we also include variation estimator  $\hat\sigma_{c,n}=\sqrt{\sum_{k=1}^n(X_{(k+1)h}-X_{kh})^2/(nh)}$, which is observed not 
consistent. 
\begin{table}
\caption{The estimators $(\hat\kappa,\hat\theta,\hat\sigma,\hat\sigma_c)$ for different values of $n$.}\label{tab:1} 
\centering   
\begin{tabular}{lllllllll}\toprule
       & \multicolumn{5}{l}{n($\times 10^3$) }\\ \cline{2-7}  
                        &2 &3&4&5&6 &8   \\ \midrule
$ \hat\kappa$&0.963&0.953&1.134&0.994&0.956&0.966\\
 $\hat\theta$ &0.992&1.001&0.996&0.998&0.989&0.997\\
 $\hat\sigma$&0.486&0.497&0.517&0.501&0.503&0.501\\
 $\hat\sigma_c$&0.431&0.443&0.451&0.444&0.446&0.444\\
\bottomrule
    \end{tabular}
\end{table}

%
%

\section*{Acknowledgement}
 Y. Hu is supported by an NSERC discovery grant and a startup fund of University of Alberta.    Y. Xi is supported by the National Natural Science Foundation of China  (Grant No.  11631004, 71532001) and the China Scholarship Council.
\bibliography{mybibfile}

\end{document}